\documentclass[a4paper,11pt]{amsart}
\usepackage{amssymb}
\usepackage[dvips]{graphicx}
\usepackage{psfrag}
\usepackage{amscd}
\newcommand{\comm}[1]{}

\def\({\left(}
\def\){\right)}

\def\raw{\rightarrow}

\def\no={\neq}
\def\sm{\setminus}

\def\C{{\mathbb C}}
\def\D{{\mathbb D}}

\def\MM{{\mathcal M}}

\def\al{\alpha}
\def\be{\beta}
\def\ga{\gamma}

\def\th{\theta}

\theoremstyle{plain}
\newtheorem{Main}{Theorem}

\newtheorem{Thm}{Theorem}[section]

\newtheorem{Lem}[Thm]{Lemma}

\newtheorem{Cor}[Thm]{Corollary}
\newtheorem*{Fact}{Fact}

\theoremstyle{remark}

\newtheorem{Def}[Thm]{Definition}

\begin{document}
\begin{center}
\end{center}
\title[The dual nest]
{The dual nest for degenerate Yoccoz puzzles}
\author{Magnus Aspenberg}

\address{Mathematisches Seminar, Christian-Albrechts Univertit\"at zu
  Kiel, Ludewig-Meyn Str.4, 24 098 Kiel, Germany}
\email{aspenberg@math.uni-kiel.de}
\thanks{The author gratefully acknowledges funding from the Research
  Training Network CODY of the European Commission}
\begin{abstract}

The Yoccoz puzzle is a fundamental tool in Holomorphic Dynamics. The
original combinatorial 
argument by Yoccoz, based on the Branner-Hubbard tableau, counts the
preimages of a non-degenerate 
annulus in the puzzle. However, in some important new applications of
the puzzle (notably, matings 
of quadratic polynomials) there is no non-degenerate annulus. We
develop a general combinatorial 
argument to handle this situation. It allows to derive corollaries,
such as the local connectedness 
of the Julia set, for suitable families of rational maps. 



\end{abstract}

\maketitle

\section{Introduction}

The Yoccoz puzzle is by now well-known in Complex Dynamics, as a way
to prove local connectivity of Julia sets for non-renormalisable (or
finitely renormalisable)
quadratic polynomials. The result was proven in 1990 by J-C. Yoccoz
(see \cite{Hub}, \cite{Mil}),
using a tableau developed by Branner and Hubbard \cite{Br-Hu}. Other
approaches has been developed also for some infinitely
renormalisable polynomials
(see e.g. \cite{ML2}, \cite{LK}, \cite{LK2} etc). The set of puzzle
pieces is a dynamical partition of a neighbourhood of the Julia set
for some function, which in the original setting was a quadratic
polynomial of the form $f_c(z)=z^2+c$, $c \in \C$. Each puzzle piece is simply connected and contains a connected part of the Julia set of $f_c$. Hence, by definition, if the puzzle pieces shrink to points, then the Julia set must be locally connected at that point. Two things must be satisfied for this to happen: Firstly, the puzzle must have a good combinatorics. This property is inherited from the fact that $f_c$ is non-renormalisable. Secondly, nested puzzle pieces must not touch each other, i.e. there must be a {\em non-degenerate} annulus between two puzzle pieces in the Branner-Hubbard tableau.

In the case of quadratic polynomials, this non-degeneracy is never a
problem, because in this case the puzzle pieces are formed by equipotentials and
external rays which behave nicely. However, the situation with degenerate
annuli appears naturally in more general 
applications of the Yoccoz puzzle, in particular for rational
functions, the reason being that the corresponding equipotentials and
external rays for rational maps have a more complicated
structure than for polynomials (see e.g. \cite{AY}). 

In this paper we give a method of how to deal with this degeneracy
problem in a consistent way. We show that if one has a Yoccoz puzzle with
"good combinatorics'' but without the non-degeneracy
condition, then the puzzle pieces still
shrink to points. The idea is to look at the space between degenerate
annuli (complementary annuli) and compute their total modulus given
that the combinatorics is correct. 

\subsection{Abstract statement of result}

We first state the main result in a combinatorial and somewhat
abstract way. In Section \ref{yoccozpuzzle} we discuss the Yoccoz puzzle and from
this it will be clear where these abstract statements come from. In
Section \ref{proof} we prove Theorem \ref{mainthm} and in Section \ref{last} we
discuss some applications.

A {\em degenerate annulus} will be understood as the closed
set bounded by two
closed Jordan arcs $\ga_1$ and $\ga_2$, which touch each other in at least one
point. We assume the interior is non-empty.
A non-degenerate annulus is simply an annulus with positive modulus.
For an annulus $A$ (degenerate or non-degenerate)
we say that the inner component $in(A)$ of $A$ is the
bounded component of $A^c$ and the outer component $out(A)$ of $A$ is the
unbounded component of $A^c$. We say that a sequence of annuli $A_j$
is {\em a-nested} if $in(A_{j+1}) \cup A_{j+1} \subset in(A_j)$
for all $j \geq 0$. Of
course the $A_j$ has to be disjoint if they are a-nested. We say that
an annulus $A$ {\em surrounds} a set $E$ if $E \subset in(A)$. Given an analytic
function $f$, a sequence of annuli
$A_j$ satisfies the Markov property if $f^n(A_j) \cap A_i = \emptyset$
unless $f^n(A_j)=A_i$, for all $i,j,n \geq 0$.

The setup is an a-nested sequence of annuli mapped onto each other in a certain
unbranched way. In the following theorem, we say that $f: A \raw B$ is a covering map
between the closed annuli $A$ and $B$, if
$f$ is a covering map on some neighbourhood of $A$ onto
some neighbourhood of $B$.
Here is the abstract formulation of the main result.

\begin{Main} \label{mainthm}
Suppose $A_j$ is an a-nested sequence of (not necessarily
non-degenerate) annuli $A_j$ satisfying the Markov property and
surrounding a critical point $z_0$ of order $2$ of an analytic
function $f$, having the following properties:
For each $A_j$, $j > 0$, there is some $A_{j'}$ and $n_j > 0$ such that
$$f^{n_j}: A_{j} \raw A_{j'}$$
is an unbranched covering of degree $2$.
Conversely,  each $A_j$ has (at least) $2$ preimages $A_{j_1}, A_{j_2}$, with $j_1, j_2 > j$, such that
$$f^{n_{j_i}} : A_{j_i} \raw A_j, \quad i=1,2 $$
are both unbranched coverings of degree $2$.

Finally, suppose that the {\em complementary annuli} $\al_j =
in(A_{j}) \cap out(A_{j+2})$, $j=0,1, \ldots$, are all
non-degenerate. Then there is a sequence of non-overlappling $\al_{j_k}$,
$k = 0,1,\ldots$, such that
\[
\sum_k \mod (\al_{j_k}) = \infty.
\]

\end{Main}

We call the set of complementary annuli a {\em dual
  nest} of annuli around the critical point. The reason for taking
the space between every second annulus in the $a$-nested sequence
$A_j$ instead of just $in(A_j) \cap out(A_{j+1})$ is that the latter
complementary annuli might be degenerate. However, the former annuli
will always be non-degenerate as we shall see.


\subsection*{Acknowledgements} I thank Carsten Petersen for
discussions which led to this paper.
I am thankful to M. Yampolsky for discussing this problem and giving
helpful remarks.
The paper was written at Mathematisches Seminar,
Christian-Albrechts Universit\"at zu Kiel. I gratefully acknowledge
the hospitality of the department.

\section{The Yoccoz puzzle} \label{yoccozpuzzle}

Let us in this section recapitulate the idea of Yoccoz famous result
of proving that the Julia set of a non-renormalisable quadratic
polynomial $f_c(z)=z^2+c$ is locally conneced. We will follow the
exposition in \cite{Mil}.
The setup is a puzzle partition by B. Branner and J. Hubbard originally made for polynomials of degree $3$. By utilising the combinatorics of the puzzle for a non-renormalisable polynomial, Yoccoz could then prove that these puzzles pieces shrink to points, thereby proving that the Julia set is locally connected.

So let $f_c(z)=z^2+c$ be non-renormalisable. It has $2$ fixed points, of which at least one must be repelling. Let $\Psi: \hat{\C} \sm \D \raw \hat{\C}\sm J(f)$ be the B\"ottcher coordinates around $\infty$. An external ray with angle $\th$ is defined by
\[
\ga(\th) = \{ z=re^{2\pi i \th} : \Psi^{-1}(re^{2 \pi i \th}), r > 1\}.
\]
An external ray with angle $\th$ lands, i.e. the limit $\lim_{r \raw
  1} \Psi^{-1}(re^{2\pi i \th})$ exists, when $\th$ is a rational
number. Evidently, landing points of such rays have to be periodic
points. Conversely, every periodic point is a landing point of
finitely many external rays with rational angles (although this is
non-trivial to prove, see e.g. \cite{Milnor-book}).

The fixed point on which an external ray with angle $\th \neq 0$ lands at, is called the $\al$-fixed point. The other fixed point, where the external ray with angle $0$ lands, is called the $\be$-fixed point. An equipotential of level $r > 1$ is defined by
\[
E(r) = \{ z = r e^{2 \pi i \th} : 0 \leq \th \leq 2 \pi \}.
\]
Assume now that the rays $\ga(\th_j)$ land at the $\al$-fixed point. Take some $r_0 > 1$. Then these rays together with the equipotential of level $r_0$ form a puzzle (of depth zero) where each puzzle piece is a bounded component of the complement of $E(r_0) \cup \cup_j \ga(\th_j)$.
Since the set $\ga(\th_j)$ of external rays are forward invariant, we
can pullback these puzzle pieces and thereby obtain the puzzle of
depth one, consiting of the puzzles pieces so that each puzzle piece
is a preimage of a puzzle piece of depth zero. If $P_d(z)$ is the
puzzle piece of level $d$ containing $z$, then it follows by
definition that $f_c(P_d(z))=P_{d-1}(f_c(z))$. We can also number the
puzzle pieces of depth $d$ and write $P_d^j$, $j=0,\ldots,n-1$ where
$n$ is the number of puzzles pieces of depth $d$. We usually denote by
$P_d^0=P_d(0)$ the critical puzzle piece containing the critical point $z=0$. Because of the forward invariance we have the following immediate Markov property:
\begin{Lem}
Given two puzzle pieces $P$ and $Q$ either $P \cap Q = \emptyset$ or
one is contained in the other.
\end{Lem}

Note that for each puzzle piece $P$ we have that $P \cap J(f_c)$ is connected. Hence if the puzzle pieces around a given point $z$ shrink to a single point, namely $z$, then the Julia set is locally connected at $z$.

To prove this we will study the annuli around the critical point $z=0$:
\[
A_d(z) = P_d(z) \sm \overline{P_{d-1}(z)}.
\]
Given an (open) annulus $A=P_d(0) \sm \overline P_{d'}(0)$ (where $d <
d'$) surrounding the critical point $z=0$ (i.e. the inner puzzle piece contains $z=0$) suppose it is subdivided into two (open) annuli $B$ and $C$ which also surround the critical point, i.e. $\overline{A} = \overline {B \cup C}$. Then we have the following Gr\"otzsch inequality:
\[
\mod(A) \geq \mod(B) + \mod(C).
\]

Hence if $\sum_d A_d(0) = \infty$ we have that $\cap_d P_d(0) = \{ 0 \}$.

Although Yoccoz's result is stated for quadratic polynomials, the ideas are
mostly combinatorial and since we will aim for more general applications, we formulate his theorem in a more combinatorial way. But first we have to define the Branner-Hubbard tableau.
Given a point $z \in J(f_c)$, not in $\cup_{j \geq 0} f_c^{-j}(\al)$,
where $\al$ is the $\al$-fixed
point, we consider its orbit under $f_c$:
\[
z \raw z_1 \raw z_2 \raw \ldots,
\]
where $f_c(z_j)=z_{j+1}$. For a given depth $d$, we say that $z_j \in P_d(z_j)$ is critical, semi-critical or off-critical whenever respectively

\begin{itemize}

\item $z_j \in P_d(0)$, (critical)

\item $z_j \in P_{d-1}(0) \sm P_d(0)$, (semi-critical)

\item $z_j \notin P_{d-1}(0)$. (off-critial)

\end{itemize}

Now note that if $z_j \in P_d(z)$ then the map $f_c: P_d(z_j) \raw P_{d-1}(z_{j+1})$ is either a double covering if $z_j$ is critical and univalent if $z_j$ is semi-critical or off-critical.

Moreover, in the critical and off-critical case, the map $f_c: A_d(z_j) \raw A_{d-1}(z_{j+1})$ is a covering map. In the semi-critical case it is not a covering map, since $P_d(z_{j+1})$ has two preimages in $P_{d-1}(z_j)$, but we have the following relationships of the modulus of these annuli:

\begin{Lem}
With notations as above, we have
\begin{align}
\mod(A_{d-1}(z_j)) &= 2 \mod (A_d(z_j)) \text{ if $z_j$ is critical}, \\
\mod(A_{d-1}(z_j)) &= \mod(A_d(z_j)) \text{ if $z_j$ is off-critical and} \\
\mod(A_{d-1}(z_j)) &< 2 \mod(A_d(z_j)) \text{ if $z_j$ is semi-critical. }
\end{align}
\end{Lem}

In the critial case above, we say that $A_{d}(z_j)$ is a {\em child} to $A_{d-1}(z_{j+1})$, i.e. when $f_c: A_d(z_j) \raw A_{d-1}(z_{j+1})$ is a double covering. The child is {\em excellent} if it has two other children.

Now we are ready to define the tableau.
\begin{Def}
The tableau is associated to a starting point $z=z_0$ and a
$2$-dimensional array of the non-negative numbers $\mu_{ij}=\mod
(A_i(z_j))$ (with the obvious notion of being critical, semi-critical
and off-critical). We mark each entry in the tableau with a critical
if $\mu_{ij}$ is critical, semi-critical
if $\mu_{ij}$ is semi-critical and off-critical
if $\mu_{ij}$ is off-critical.
\end{Def}
A movement in the north-east direction in the tableau, from $\mu_{ij}$
to $\mu_{i-1,j+1}$, represents the action of the map $f_c$ on the
annulus $A_i(z_j)$. It follows from the definition that each column in
the tableau starts with critial entries, then at some point comes a
semi-critical entry, and below this only off-critical entries exist.
These entries are usually depicted as single lines (critial marking),
double lines (semi-critial marking), and no lines (off-critial
marking).

We will from now on only consider the critical tableau, i.e. when $z_0=c_0=0$. Write $c_0=0, c_1=c, c_n=f_c(c_{n-1})$.
We say that the critical tableau is {\em recurrent} if
\[
\sup \{ d : \mu_{dk} \text{ is critical }  \} =  \infty.
\]
We say that the critical tableau is {\em periodic} if some $k$th column, $k > 0$, is entirely critical.

\begin{Thm}(Yoccoz)
Assume that the critial tableau is recurrent but not periodic and that there exists some non-degenerate annulus $A_d(0)$ such that $\mod (A_d(0)) > 0$. Then \[
\sum_d A_d(0) = \infty.
\]
\end{Thm}

This result depends essentially on the following two facts.
\begin{Lem} \label{children}
Let $\mu_{ij}$ be a critical tableau which is recurrent but not periodic.

\begin{enumerate}

\item Assume there is a child $\mu_{d0}=A_d(0)$ which is excellent. Then all children to $A_d(0)$ are excellent.

\item There exists at least one child $A_d(0)$ which is excellent in the tableau.
\end{enumerate}

\end{Lem}

If $f_c$ is non-renormalisable, then
some $A_d(0)$ is an excellent child and
this child's all descendants are all
excellent. If $A_{d'}(0)$ is a child to $A_d(0)$ we have $\mod
(A_d(0)) = 2\mod (A_{d'})$, because $f_c: A_{d'}(0) \raw A_d(0)$ is a
double unbranched covering. Since every excellent child has at least $2$
children, and the modulus of each of these children is half of their parents,
we get that the total sum of the moduli is
\[
\sum_{A_{d'}(0) \text{ is a descendant to }A_d(0)} A_{d'}(0) \geq \sum_{k \geq 0} 2^k
\frac{1}{2^k} \mod (A_d(0)) = \infty,
\]
given that the top child $A_d(0)$ is
non-degenerate, i.e. $\mod (A_d(0)) > 0$. The existence of such a
non-degenerate child is
automatic as soon as some $A_n(0)$ is non-degenerate.

Hence Yoccoz's result follows from \ref{children} given the
non-degeneracy condition.
But what happens if $A_d(0)$ and all its children are degenerate? The
main result of the paper is that the puzzle pieces still shrinks to
points.

\section{The complementary annuli and Proof of
  Theorem \ref{mainthm}.} \label{proof}

From the previous section we have seen that the assumptions in
Theorem \ref{mainthm} is natural and comes from the construction of
the Yoccoz puzzle. In particular the sequence of annuli $A_j$ is
simply the set of descendants of an excellent child in the puzzle.
Now let us prove Theorem \ref{mainthm}.

Assume $A=A_0$ is a degenerate critical annulus
and that $A_j$, $j=0,\ldots$ satisfies the assumptions of the
theorem.
Hence every child $A_j$ has at least two
children, each mapped onto $A_j$ as an unbranched convering of degree
$2$.
Let us relabel these annuli. They form a tree of descendants $A_{i,j}$ starting from $A=A_{0,1}$ so that, for fixed $i > 0$, $A_{i,j}$ are the descendants of generation $i$. Generation $i$ means that $f^i(A_{i,j})=A_0$ and that $f^k: A_{i,j} \raw A_0$ is a $2^i$ degree unbranched covering. Moreover, since every $A_{i,j}$ is excellent there are at least $2^i$ annuli of generation $i$.


A {\em complementary annulus} $\al_j$ is defined by the annulus
bounded by $A_j$ and $A_{j+2}$. The annulus $A_{j+1}$ is called the {\em
  middle annulus} (between $A_j$ and $A_{j+2}$) of $\al_j$. Of course
these annuli
overlap unless we take every second annulus. We will deal with this later.

In this proof we always assume that a complementary
annulus is {\em non-degenerate}. In the applications it turns out that
they are. We want to see
what kind of relation there is between $A_j$ and the $\al_j$.




Take some complementary $\al$ bounded by the degenerate annuli
$P$ and $Q$, where $P$ surrounds $Q$. Note that we assume that exactly
one annulus $R$ lies between $P$ and $Q$.
Now $Q$ has a child, say $Q_1$, so that $Q_1$ maps onto $Q$ as a $2$ degree unbranched covering. We want to pull back $P$ along the same branch (if possible) as $Q$ back to some $P_j$ surrounding $Q_1$.

In the first steps $\al$ (between $P$ and $Q$) is pulled back as a
one-to-one map until some preimage $P_1$ of $P$ under $f^{k}$
surrounds the critical point $z_0$.
This means by definition that this preimage $P_1$ is a child to
$P$. If moreover $Q_1$, being the preimage of $Q$ under $f^{k}$
surrounded by $P_1$, also surrounds the critical point we are done and
have found $P_1$ surrounding $Q_1$ both being children of $P$ and $Q$
respectively. Clearly, if exactly one degenerate annulus lies between
$P$ and $Q$ then exactly one degenerate annulus lies between $P_1$ and $Q_1$.

The second (and most probable) case is however that, whereas $P_1$ surrounds the critical point, $Q_1$ does not surround the critical point. Hence we are in a semi-critical situation, so the pullback $f^{-k}(\al)$ is not an annulus. However, if we consider the annulus $\be_1$ between $P_1$ and $Q_1$, this annulus has modulus at least $1/2$ of the modulus of $\al$, (by standard inspection from semi-critical annuli). Continuing pulling back $\be_1$, we again sooner or less reach the same situation:
Some pullback $P_2$ of $P_1$ under $f^{k_1}$ surrounds the critical
point. If again the preimage $Q_2$ (being a preimage of $Q_1$ under
$f^{k_1}$) surrounded by $P_2$ also surrounds the critical point we
are done and have found two descendants $P_2$ and $Q_2$ to $P$ and $Q$
respectively. However, note that, where as $Q_2$ is a child to $Q$, we
have that $P_2$ is a child of $P_1$ and $P_1$ is a child of
$P$. ($Q_1$ is not a child of $Q$ since $Q_1$ was assumed not to
surround the critical point).

Continuing in this way we find two descendants $P_m$ and $Q_m$ such that
$$f^{k+k_1+\ldots+k_{m-1}}: P_m \raw P$$
as a $2^m$ degree unbranched covering and
$$f^{k+k_1+\ldots+k_{m-1}}: Q_m \raw Q$$
as a $2$ degree unbranched covering.

In the same way, the middle annulus $R$ between $P$ and $Q$ is pulled
back to some $R_m$ between $P_m$ and $Q_m$ and
$$f^{k+k_1+\ldots+k_{m-1}}: R_m \raw R, $$
as a $D$ degree unbrached covering where $2 \leq D \leq 2^m$.

Here $Q_m$ is a child to $Q$, whereas every $P_{j+1}$ is a child to
$P_j$, $j=0, \ldots, m-1$. The annulus $R_{j+1}$ is a child to $R_j$ only if $R_{j+1}$ surrounds the critical point.

We call the annulus bounded by
$P_m$ and $Q_m$ an {\em offspring} to $\al$ and $\al$ the {\em
  ancestor } to the annulus between $P_m$ and $Q_m$.
Hence every offspring has modulus at least $2^{-m}$ times the modulus
of its ancestor $\al$, where $m$ is defined above.

Obviously, if there is a degenerate annulus $A$ between $P_m$ and
$Q_m$, we can map this annulus forward;
$f^{k+k_1+\ldots+k_{m-1}}(A)$ will in that case be some degenerate
annulus between $P$ and $Q$.

Conversely, let $P_m$ and $Q_m$ be given degenerate annuli bounding
the complementary annulus $\al_1$ and assume that exactly one degenerate
annulus $R_m$ lies between $P_m$ and $Q_m$.
If $Q_m$ has generation more than $1$ then the parent $Q$ would have
generation at most $1$. On the other hand, the parent $P$ to $P_1$,
which in turn is parent to $P_2$ and so on down to $P_m$, might have negative
generation, meaning that $P$ is actually a parent to $A_0$. In this
case, $A_0$ would lie between $P$ and $R$. But in this case there has
to be some preimage $A_j'$ of $A_0$ laying between $P_m$ and $Q_m$. If
also $R$ has negative generation, then there is some preimage to $A_0$
different from $R_m$ between $Q_m$ and $R_m$, a contradiction. Hence
$R$ has generation at most zero.

If the generation of $R$ is greater
than zero, then the generation of $P$ has to be zero, i.e. $P=A_0$,
because otherwise we could pullback $A_0$ to some degenerate annulus
between $R_m$ and $P_m$.

From this we see conclude:
\begin{Fact}
If there is exactly one degenerate annulus between $Q_m$ and $P_m$
then there is exactly one degenerate annulus between $Q$
and $P$, if the generation of $R$ is at least $1$.
\end{Fact}

\begin{Def}
Given a complementary annulus $\al$ bounded by $P_m$ and $Q_m$, with
middle annulus $R_m$, we say that its {\em intermediate generation} is
equal to the generation of $R$, where $R = f^{k+k_1+ \ldots +
  k_{m-1}}(R_m)$ and $f^{k+k_1+\ldots +k_{m-1}}$ is the function
described above.
\end{Def}

From the discussion we conclude:
\begin{Lem}
Every complementary annulus $\al$ with intermediate generation at
least $1$ has some unique ancestor $\be$.
\end{Lem}
\begin{proof}
Clearly, the annulus $P=f^{k+k_1+\ldots + k_{m-1}}(P_m)$ must have
generation at most zero. That means that no denegerate annulus can
exist between $P$ and $R$ or between $R$ and $Q$ because this would
then have a preimage $R'$ between $Q_m$ and $P_m$ different from $R$,
a contradiction. Hence the annulus $\be$ bounded by $P$
and $Q$ is indeed a complementary annulus.
\end{proof}

\begin{Def}
Given a complementary annulus $\al$ bounded by the outer degenerate
annulus $A_{m,*}$ and the inner degenerate annulus $A_{n,*}$,
we say that the {\em outer generation} to $\al$ is equal to $m$ and the {\em inner generation} to $\al$ is $n$. We write $\al=\al_{n,*}^m$, where $*$ means an index, since there might be many $\al$ with the same $m$ and $n$.

\end{Def}

We have proved the following.

\begin{Lem} \label{onestep}
For every complementary annulus $\al=\al_{n,*}^m$ with intermediate
generation at least $1$ and with ancestor $\al_{n-1,*}^{m_1}$ we have
\[
\mod(\al_{n,*}^m) \geq 2^{m_1-m} \mod(\al_{n-1,*}^{m_1}).
\]
\end{Lem}

\begin{Cor} \label{manysteps}
For every complementray annulus $\al_{n,*}^m$, $n > 1$ with
intermediate generation at least $1$,
there is some grand ancestor $\al_{N,*}^{m_{n-1}}$ such that
\[
\mod (\al_{n,*}^m) \geq 2^{m_{n-N}-m} \mod (\al_{N,*}^{m_{n-N}}).
\]
\end{Cor}

Since the number of degenerate annuli of generation $m$ is at least $2^m$ we have that the number of complementary annuli of outer generation $m$ is at least $2^m$.
If the complementary annuli are non-degenerate, there is some $M_0 >
0$, such that $\mod(\al_{N,*}^m) \geq M_0$ for all grand ancestors $\al_{N,*}^m$.


When we consider the complementary annuli, we want to sum every second
moduli (since otherwise they overlap). Clearly, the set of complementary annuli which have intermediate
generation at most $0$ are finite.  Therfore, we can fix a generation
$m_0 > 1$ such that all $\al_{n,*}^{m}$, $m \geq m_0$, have intermediate generation
at least $1$. So every such annulus has some grand ancestor.

Now consider the set of the annuli $\al_{n,*}^{m_0}$. Either half of
them will have the property
that the inner annulus is $A_j$ where $j$ is even or at least half of them
will have an inner annulus being $A_j$ where $j$ is odd. Suppose that
the first case occurs. In this collection of complementary annuli all the inner generations are even so they do not overlap. We get
\[
\sum_{n,*} \mod (\al_{n,*}^{m_0}) \geq 2^{m_0-1} 2^{m_{n-N}-m_0} \mod (\al_{N,*}^{m_{n-N}}) \geq M_0/2,
\]
where the sum runs over the annuli $\al_{n,*}^{m_0}$, where all of
them have inner generation even. Of course the same statement holds in
the odd case.

Going sufficiently deep in the nest we pick another $m_1 > m_0$ such
that no $\al_{n,*}^{m_1}$ intersects any $\al_{n,*}^{m_0}$. Again, the sum of
these annuli which do not overlap becomes at least $M_0/2$.
Continuing in this manner we get a sequence
$m_0 \leq m_1 \leq m_2 \ldots$ of outer generations such that the sum
of non-overlaping complementary annuli of generations $m_0, m_1,
\ldots$ becomes at least $$M_0/2 + M_0/2 + \ldots
= \infty,$$
and Theorem \ref{mainthm} follows.

\section{The non-degeneracy of the complementary annuli} \label{last}

So far we have seen that degenerate annuli in the Yoccoz puzzle is not
an obstacle to prove that puzzle pieces shrink to points, under the
condition that the complementary annuli are non-degenerate. In the applications, the combinatorics comes from the Mandelbrot set $\MM$ and it turns out that the complementary annuli are automaticly non-degenerate.

Given a non-renormalisable quadratic polynomial $f_c(z)=z^2+c$, the schematic
picture of the first two levels of the Yocccoz puzzle looks like in
the figure above, (here we have chosen
$c$ not from the $1/2$-limb of the Mandelbrot set).

\begin{figure} \label{puzzle1}
\psfrag{a}[][][0.9]{$\al$}
\psfrag{b}[][][0.9]{$-\al$}
\psfrag{P00}[][][0.9]{$P_0^0$}
\psfrag{P11}[][][0.9]{$P_1^1$}
\centerline{\includegraphics[width=1.0\textwidth]{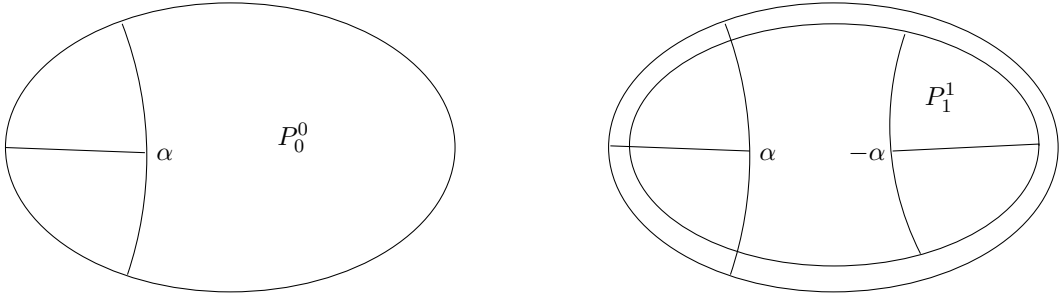}}
\caption{The Yoccoz puzzles of depths $0$ and $1$, with $q=3$ external rays
  landing at $\al$. The ellipse to the left is an equipotential and the inner ellipse to the right is its preimage.}
\end{figure}

Generally, puzzle pieces containing the critical point are mapped onto
each other
under iterations of $f^q$, where $q \geq 2$,
i.e. at least $2$ iterates of $f$ is
required for some critical puzzle piece $P_d^0$ to be mapped onto
another critical puzzle piece $P_{d'}^0$ (in the figure $q = 3$). Hence, between $P_d^0$ and
$P_{d'}^0$ there are at least one puzzle pieces containing $P_{d'}^0$
and contained in $P_d^0$. If now, the two top puzzle pieces $P_1^1$
 and
$P_0^0$ touch each other at some set $E$, the set $E$ will be pulled
back and produce degeneracies at all depths between preimages of
$P_1^1$ and $P_0^0$. However, these touching points cannot coincide
and we have the following:

\begin{Fact}
For any depth $d \geq 0$ and nested puzzle pieces $P_d \supset P_{d+1}
\supset P_{d+2}$ we have that $P_{d+2} \Subset P_d$, i.e. there are
no touching points between $P_d$ and $P_{d+2}$.
\end{Fact}

So if we go two
levels down in the nest we create a non-degenerate annulus. Since we
always go at least $2$ leves down between consecutive critical puzzle
pieces, the complementary annuli $\al_d$ being the space between
$A_d(0)$ and its pullback to one if its grand children $A_{d'}(0)$, $d' \geq
d+2q \geq d +4$, have to be non-degenerate, i.e. $\mod (\al_d) > 0$.

Note that if $c$ belonged to the $1/2$-limb of $\MM$, then the
annulus between $A_j$ and $A_{j+1}$ could be {\em degenerate}
also. This is the reason why we consider the complementary annuli
being the space between every second annulus in the $a$-nested sequence $A_j$, $j=0,\ldots$.

Hence for all non-renormalisable combinatorics from the Mandelbrot set Theorem \ref{mainthm}
works. In particular, the puzzle pieces in the Yoccoz puzzle $P_d(z_0) = in (A_{d-1}(z_0))$ containing the critical point $z_0$ must shrink to a single point:
$$\bigcap_{j\geq 0} in(A_j) = \{z_0 \}. $$


\comm{
\subsection{Combinatorics from the $1/2$-limb of $\MM$.}
If $c$ belongs to the $1/2$-limb of $\MM$,
then the critical puzzle pieces are mapped onto each other under
iterations of $f^q$, where $q \geq 2$. Hence the descendants of some
(degenerate) $A_d(0)$ may touch each other and the simple complementary
annuli may be degenerate. However, in this case we look at
every second annulus in the descendants $A_j$ of $A_0=A_d(0)$,
i.e. the enlarged annuli. By the second statement in Theorem
\ref{mainthm}, we still have infinite sum of their moduli.
}


\bibliographystyle{plain}
\bibliography{ref}

\end{document}